\documentclass[preprint,12pt]{elsarticle}

\usepackage{amssymb}

\begin{document}

\begin{frontmatter}



\title{Attractiveness of Invariant Manifolds\tnoteref{l1}}
 \tnotetext[l1]{This work was supported by the NNSF of China
under the Grant No. 10702065.}



\author[aa1]{Lijun Pei \tnoteref{l2}}

\address[aa1]{Department of Mathematics,  Zhengzhou University, 450001
Zhengzhou, Henan, China}\tnotetext[l2]{Corresponding Author. Tel:
+86 371
  67783167.\\
 Email Address: peilijun@zzu.edu.cn, lijunpei@yahoo.com.cn.
\newline }


\begin{abstract}
In this paper an operable, universal and simple theory on the attractiveness of the
invariant manifolds is first obtained. It is motivated by the Lyapunov direct method. It means that for any point $\overrightarrow{x}$ in the invariant manifold
$M$, $n(\overrightarrow{x})$ is the normal passing by $\overrightarrow{x}$, and $\forall \overrightarrow{x^{'}} \in  n(\overrightarrow{x})$, if the tangent $f(\overrightarrow{x^{'}})$ of the orbits of the
dynamical system intersects at  obtuse (sharp) angle with the normal  $n(\overrightarrow{x})$, or
the inner product of
the normal vector $\overrightarrow{n}(\overrightarrow{x})$ and tangent vector $\overrightarrow{f}(\overrightarrow{x^{'}})$ is
negative (positive), i.e.,
$\overrightarrow{f}(\overrightarrow{x^{'}}).
\overrightarrow{n}(\overrightarrow{x}) < (>)0$, then the invariant manifold $M$ is attractive (repulsive). Some illustrative
examples of the invariant manifolds, such as equilibria, periodic solution, stable and unstable manifolds,
other invariant manifold are presented to support our result.
\end{abstract}

\begin{keyword}


attractiveness;\ invariant manifold;\ equilibria;\  periodic solutions;\ stable and unstable
manifolds
\end{keyword}

\end{frontmatter}


\section{Introduction}

\vspace*{0.2cm}

The theory of invariant manifolds (shorten for IMs) is very important to the
reduction of the higher-dimensional and complex systems, synchronization of
the coupled and complex chaotic systems. The existence, uniqueness and multivaluedness of the IMs can
be solved by the theory of Partial Differential Equations (shorten for PDEs) since the differential
equations governing the IMs are in fact the first-ordered
quasi-linear PDEs. So the existence and number of IMs are equivalent to
 those of the analytic solutions of  the first ordered quasi-linear PDEs.
 The existence and uniqueness of the IM can be determined by Cauhy-Kawalewskaja Theorem $^{\mathbf{[1]}}$.
  And  several IMs appear if the above theorem doesn't hold. Since the general analytic
  solutions of  the first ordered quasi-linear PDEs can't be solved explicitly, the expression of IM
  can be approximated only by the numerical method. The another important question of the IMs is
 their attractiveness.

 The attractiveness of the IMs is a difficult and unsolved question.
 Some authors present the different methods to solve it.
 Fenichel $^{\mathbf{[2]}}$ obtained the sufficient conditions for persistence of a
 diffeomorphic IM under perturbation of the flow in terms of generalized Lyapunov type numbers, and the smoothness of the perturbed manifold.
 The concept of normally hyperbolicity has been
 proposed in $\mathbf{[2]}-\mathbf{[5]}$, i.e.,
 the contraction of the flow is $r \geq 1$ times exponentially stronger in the direction normal to the manifold than within the manifold.
 Josi$\acute{c}$ introduced a modification of Fenichel theory which applied to chaotic synchronization,
  proposed a necessary and sufficient condition for such persistence was $r$ normally hyperbolicity,
  and discussed the Lyapunov-exponent-like quantities used to determine the transverse stability of synchronization manifolds $^{\mathbf{[6]}}$.
 But these methods
 are too abstract to be employed. The other gave the operable
 but not general method to attack it, for example,  Gorban  put forward
 the approach that the stability of  the equilibria of the invariance equation, which
 correspond  to   the slow IMs, is equivalent to that or attractiveness  of the
 corresponding slow IMs, thus the  slow IMs' stability
  can be obtained
 by the corresponding stability $^{\mathbf{[7]}}$. There is a generalization that the
 IMs' stability or attractiveness is also
 equivalent to stability of the
 corresponding   equilibrium of the invariance equation.
  But if the IM is not also the
 equilibrium of the invariance equation, such as the global periodic solutions,  the attractiveness can't be
 solved in this way. Thus an operable, simple  and general method is looked forward to
 appearing to consider the attractiveness of the  IMs.

In this paper a new method of the attractiveness of the  IMs is first proposed. The attractiveness of the  IMs is  equivalent to
that the nearby orbits will intersect continuously inward with the normals of
IM, or the tangent of the orbits passing the normals intersect
at  obtuse angle. Obviously the latter suggests the inner product
of the  tangent vectors of the orbits and the normal vectors of
IM is negative. This idea is motivated from the attractiveness of the equilibria by the theory of
Lyapunov direct method $^{\mathbf{[8]}}$. Then some examples are presented to verify
this idea.

 The structure of this paper is the following: the theory of the attractiveness of the IMs
 of the two-dimensional dynamical systems is derived in Section 2; then some illustrative examples
 are presented in Section 3;
 at last the conclusion and discussion are given.

\vspace{0.2in}
 \section{Attractiveness of the IMs
 of the two-dimensional dynamical systems}
\vspace*{0.3cm}

There are several equivalent arguments on the attractiveness of the IMs
 of the two-dimensional dynamical systems. And  an operable,
 simple  and general method can be achieved by them.
 If the IM is attractive, then the orbits in some
neighbor field will be attracted to it and close to it, i.e.,
the distance of the nearby orbits  to the IM will tend to zero as time
tends to the positive infinity. The latter is  equivalent to
that the nearby orbits will intersect continuously inward with the normals of
IM, or the tangent of the orbits passing the normals intersect
at  obtuse angle. Obviously the latter suggests the inner product
of the normal vectors of
IM and the  tangent vectors of the orbits passing the normal  is negative. This idea is motivated from the attractiveness of the equilibria by the theory of
Lyapunov direct method $^{\mathbf{[8]}}$.
It means that if the gradient $grad V$ of the normal along the equipotential
surfaces $V \equiv c$ of the $V$ function cross the tangent of the
orbit, then the orbit won't traverse outward the equipotential
surfaces $V \equiv c$, thus the equilibria are  asymptotically
stable.

$\bf{Remark 1}$.  For the attractiveness of the IMs, it is required that
the separation angle of the tangent of the orbit passing $\bf{every}$ point on the normal
of  $\bf{every}$ point on the IM and the normal is only obtuse angle but not  other angles including the right angle, i.e.,
their inner product is negative but not positive or 0.

$\bf{Remark 2}$. Here the separation angle of the tangent and the  gradient $grad V$ of the IM
is not considered, since  it judges if the orbit cross the IM. By the  uniqueness of the
solution of the ODEs,  the orbit is impossible to cross the IM. Thus we must consider
if the the orbit cross inward the normal and close to the IM.

$\bf{Remark 3}$. If the IM is the closed curve, such as the periodic solution or
limit cycle, then not only the inward attractiveness but also the  outward attractiveness must be both considered.

Then the theory for the attractiveness of the IMs is obtained.

$\bf{Theorem  1}$ Assuming the manifold $M$ is the IM of the two-dimensional dynamical system
\begin{eqnarray}
\dot{\overrightarrow{x}}=\overrightarrow{f}(\overrightarrow{x})  \label{b1}
\end{eqnarray}
,
for $\forall \overrightarrow{x} \in M$, $n(\overrightarrow{x})$ is the normal passing by $\overrightarrow{x}$, and $\forall \overrightarrow{x^{'}} \in  n(\overrightarrow{x})$,  $f(\overrightarrow{x^{'}})$ is the tangent
 of the orbits of the
dynamical system passing $\overrightarrow{x^{'}}$, then the attractiveness and repulsiveness are obtained,
\begin{enumerate}
  \item  if the tangent $f(\overrightarrow{x^{'}})$
  intersects at  obtuse angle with the normal  $n(\overrightarrow{x})$, or the inner product of
the normal vector $\overrightarrow{n}(\overrightarrow{x})$ and tangent vector $\overrightarrow{f}(\overrightarrow{x^{'}})$ is negative,   i.e., $\overrightarrow{f}(\overrightarrow{x^{'}}).
\overrightarrow{n}(\overrightarrow{x}) < 0$, then the IM $M$ is attractive;
  \item otherwise, if the tangent $f(\overrightarrow{x^{'}})$
 intersects at  sharp  angle with the normal  $n(\overrightarrow{x})$, or the inner product of
the normal vector $\overrightarrow{n}(\overrightarrow{x})$ and tangent vector $\overrightarrow{f}(\overrightarrow{x^{'}})$ is positive,   i.e., $\overrightarrow{f}(\overrightarrow{x^{'}}).
\overrightarrow{n}(\overrightarrow{x}) > 0$, then the IM $M$ is repulsive;
  \item if the tangent $f(\overrightarrow{x^{'}})$
  intersects at  right  angle with the normal  $n(\overrightarrow{x})$, or the inner product of
the normal vector $\overrightarrow{n}(\overrightarrow{x})$ and tangent vector $\overrightarrow{f}(\overrightarrow{x^{'}})$ is zero,   i.e., $\overrightarrow{f}(\overrightarrow{x^{'}}).
\overrightarrow{n}(\overrightarrow{x}) = 0$, is neither attractive nor repulsive.
\end{enumerate}

Then the attractiveness of the IMs in the two-dimensional dynamical systems, such
as that of  the equilibria, their stable and unstable manifolds, the periodic solutions (i.e., global solutions)
and other IMs, is considered by Theorem  1. It displays the correctness,
universality and simpleness of Theorem 1.

\section {Examples of attractiveness of different IMs}

\subsection {Attractiveness of equilibria}

The equilibria can be assumed to be zero without the loss of generality. Their
normal is  $y= k x, \forall k \in \mathbf{R}$ since their tangent
is the point $(0, 0)$ and they are orthogonal. The attractiveness of all kinds of
simple equilibria, i.e., the focus, node, saddle and center,  will be considered  by
Theorem  1 in this section.

 \subsubsection {Focus}
 The considered system is
 \begin{eqnarray}
 \left\{\begin{array}{l}
\dot{x}=-x + y,\\
\dot{y}=-x - y,
\end{array}\right.\label{b2}
\end{eqnarray}
$(0,0)$ is the stable focus of system $(\ref{b2})$ and attractive.
$(0,0)$ is also it's IM. The normal of the focus $(0,0)$
is $y= k x, \forall k \in \mathbf{R}$. Let $(x^{'}, y^{'})$ is
any point in the normal $y= k x$, i.e., $y^{'}= k x^{'}$. The tangent
vector of the orbit passing $(x^{'}, y^{'})$ is $\overrightarrow{f}(\overrightarrow{x^{'}})=((k-1) x^{'},(-k-1) x^{'})^{T}$.
In the
first quadrant, the normal vector is $\overrightarrow{n}(\overrightarrow{x})=(1,k)^{T}$ and
their inner product is $\overrightarrow{n}(\overrightarrow{x}).\overrightarrow{f}(\overrightarrow{x^{'}})=- (1+k^{2}) x^{'}<0$
since $x^{'}>0$. By Theorem 1, the focus is attractive in the first quadrant. In the second quadrant,
the normal vector is $\overrightarrow{n}(\overrightarrow{x})=(-1,-k)^{T}$ and the inner product is $\overrightarrow{n}(\overrightarrow{x}).\overrightarrow{f}(\overrightarrow{x^{'}})=(1+k^{2}) x^{'}<0$
since $x^{'}<0$. In the third quadrant,
the normal vector is $\overrightarrow{n}(\overrightarrow{x})=(-1,-k)^{T}$ and the inner product is $\overrightarrow{n}(\overrightarrow{x}).\overrightarrow{f}(\overrightarrow{x^{'}})=(1+k^{2}) x^{'}<0$
since $x^{'}<0$. In the fourth quadrant,
the normal vector is $\overrightarrow{n}(\overrightarrow{x})=(1,k)^{T}$ and the inner product is $\overrightarrow{n}(\overrightarrow{x}).\overrightarrow{f}(\overrightarrow{x^{'}})=-(1+k^{2}) x^{'}<0$
since $x^{'}>0$. So the focus is attractive in the other quadrants. Thus, the focus is attractive in the all quadrants
of the plane.

The repulsiveness of the unstable focus is also considered here.
For the system
 \begin{eqnarray}
 \left\{\begin{array}{l}
\dot{x}=x + y,\\
\dot{y}=-x + y,
\end{array}\right.\label{b3}
\end{eqnarray}
$(0,0)$ is the unstable focus of system $(\ref{b3})$ and repulsive.
The normal of the focus $(0,0)$
is still $y= k x, \forall k \in \mathbf{R}$. The tangent
vector of the orbit passing $(x^{'}, y^{'})$ is $\overrightarrow{f}(\overrightarrow{x^{'}})=((k+1) x^{'},(k-1) x^{'})^{T}$.
In the
first quadrant, the normal vector is $\overrightarrow{n}(\overrightarrow{x})=(1,k)^{T}$ and
their inner product is $\overrightarrow{n}(\overrightarrow{x}).\overrightarrow{f}(\overrightarrow{x^{'}})=(1+k^{2}) x^{'}>0$
since $x^{'}>0$. By Theorem 1, the focus is repulsive in the first quadrant. Thus
the focus is repulsive in the full plane.

\subsubsection {Node}
 The considered system is
 \begin{eqnarray}
 \left\{\begin{array}{l}
\dot{x}=-x,\\
\dot{y}=-2 y,
\end{array}\right.\label{b3}
\end{eqnarray}
$(0,0)$ is the stable node of system $(\ref{b3})$ and attractive.
$(0,0)$ is also it's IM. The normal of the node $(0,0)$
is still $y= k x, \forall k \in \mathbf{R}$. The tangent
vector of the orbit passing $(x^{'}, y^{'})$ is $\overrightarrow{f}(\overrightarrow{x^{'}})=(-x^{'},-2 k x^{'})^{T}$.
In the
first quadrant, the normal vector is $\overrightarrow{n}(\overrightarrow{x})=(1,k)^{T}$ and
their inner product is $\overrightarrow{n}(\overrightarrow{x}).\overrightarrow{f}(\overrightarrow{x^{'}})=-(1+2 k^{2}) x^{'}<0$
since $x^{'}>0$.
In the second and third quadrants,
the normal vector is $\overrightarrow{n}(\overrightarrow{x})=(-1,-k)^{T}$ and the inner product is $\overrightarrow{n}(\overrightarrow{x}).\overrightarrow{f}(\overrightarrow{x^{'}})=(1+2 k^{2}) x^{'}<0$
since $x^{'}<0$. In the fourth quadrant,
the normal vector is $\overrightarrow{n}(\overrightarrow{x})=(1,k)^{T}$ and the inner product is $\overrightarrow{n}(\overrightarrow{x}).\overrightarrow{f}(\overrightarrow{x^{'}})=-(1+2 k^{2}) x^{'}<0$
since $x^{'}>0$.  Thus, the node is attractive in the full plane.

The repulsiveness of the unstable node is also considered here.
For the system
 \begin{eqnarray}
 \left\{\begin{array}{l}
\dot{x}=x,\\
\dot{y}=2y,
\end{array}\right.\label{b4}
\end{eqnarray}
$(0,0)$ is the unstable node of system $(\ref{b4})$ and repulsive.
The normal of the node $(0,0)$
is still $y= k x, \forall k \in \mathbf{R}$. The tangent
vector of the orbit passing $(x^{'}, y^{'})$ is $\overrightarrow{f}(\overrightarrow{x^{'}})=(x^{'},2 k x^{'})^{T}$.
In the
first quadrant, the normal vector is $\overrightarrow{n}(\overrightarrow{x})=(1,k)^{T}$ and
their inner product is $\overrightarrow{n}(\overrightarrow{x}).\overrightarrow{f}(\overrightarrow{x^{'}})=(1+2 k^{2}) x^{'}>0$
since $x^{'}>0$. By Theorem 1, the focus is repulsive in the first quadrant. Thus similarly
the unstable node is repulsive in the full plane.

\subsubsection {Saddle}
 The considered system is
 \begin{eqnarray}
 \left\{\begin{array}{l}
\dot{x}=x,\\
\dot{y}=-2 y,
\end{array}\right.\label{b5}
\end{eqnarray}
$(0,0)$ is the saddle of system $(\ref{b5})$ and unattractive.
 The tangent
vector of the orbit passing $(x^{'}, y^{'})$ is $\overrightarrow{f}(\overrightarrow{x^{'}})=(x^{'},-2 k x^{'})^{T}$.
In the
first quadrant, the normal vector is $\overrightarrow{n}(\overrightarrow{x})=(1,k)^{T}$ and
their inner product is $\overrightarrow{n}(\overrightarrow{x}).\overrightarrow{f}(\overrightarrow{x^{'}})=(1-2 k^{2}) x^{'}$
and is not sign definite, since $k$ is varying. Thus, the saddle is unattractive.

\subsubsection {Center}
 The considered system is
 \begin{eqnarray}
 \left\{\begin{array}{l}
\dot{x}=y,\\
\dot{y}=-x,
\end{array}\right.\label{b6}
\end{eqnarray}
$(0,0)$ is the center of system $(\ref{b6})$, stable and but unattractive.
 The tangent
vector of the orbit passing $(x^{'}, y^{'})$ is $\overrightarrow{f}(\overrightarrow{x^{'}})=(k x^{'},-x^{'})^{T}$.
In the
first quadrant, the normal vector is $\overrightarrow{n}(\overrightarrow{x})=(1,k)^{T}$ and
their inner product is $\overrightarrow{n}(\overrightarrow{x}).\overrightarrow{f}(\overrightarrow{x^{'}})=0$.
In the other quadrants, their inner product is always $\overrightarrow{n}(\overrightarrow{x}).\overrightarrow{f}(\overrightarrow{x^{'}})=0$.
 Thus, the center is neither attractive nor repulsive.  It is just stable but not asymptotically stable.

\subsection {Attractiveness of periodic solution}

For the canonical system
 \begin{eqnarray}
 \left\{\begin{array}{l}
\dot{x}=-y-x (x^{2}+y^{2}-1),\\
\dot{y}=x-x (x^{2}+y^{2}-1),
\end{array}\right.\label{b7}
\end{eqnarray}
obviously there is a stable periodic solution, or limit cycle, a global solution $x^{2}+y^{2}=1$.
It's attractiveness can't be considered by the stability of
the corresponding equilibrium of the invariance equation since it isn't the equilibrium of the invariance equation.
Let $y=f(x)$ is the IM of system $(\ref{b7})$
and the invariance equation is
 \begin{eqnarray}
 \displaystyle \frac{df(x(t))}{dt}=x-f(x)[x^{2}+f(x)^{2}-1]. \label{b9}
\end{eqnarray}
In fact the stable periodic solution  $x^{2}+y^{2}=1$ isn't the equilibrium of the invariance equation
$(\ref{b9})$ and it's attractiveness can't be considered in terms of the method in $\mathbf{[7]}$.
But it's attractiveness can be deduced by Theorem 1 in Section 2.
Since it is the closed curve, not only the inward attractiveness but also the  outward attractiveness must be considered.

Let any point $(x_{0}, y_{0})$ in the periodic solution is  in the third quadrat, where
$x_{0}, y_{0} <0$. It's normal is $y=\frac{y_{0}}{x_{0}}x$. Let any point $(x^{'}, y^{'})$ in the
normal, where $y^{'}=\frac{y_{0}}{x_{0}}x^{'}$, the inward normal  vector
is $\overrightarrow{n_{i}}(\overrightarrow{x})=(1,\frac{y_{0}}{x_{0}})^{T}$.
The tangent vector of the point $(x^{'}, y^{'})$ is
 \begin{eqnarray}
\overrightarrow{f}(\overrightarrow{x^{'}})=(-\frac{y_{0}}{x_{0}}x^{'}-x^{'}(\frac{x^{'2}}{x_{0}^{2}}-1),x^{'}-\frac{y_{0}}{x_{0}}x^{'}(\frac{x^{'2}}{x_{0}^{2}}-1))^{T}.
\label{b10}
\end{eqnarray}
The inner product of the the inward normal  vector
 $\overrightarrow{n_{i}}(\overrightarrow{x})$ and the tangent vector $\overrightarrow{f}(\overrightarrow{x^{'}})$
is
\begin{eqnarray}
\overrightarrow{n_{i}}(\overrightarrow{x}).\overrightarrow{f}(\overrightarrow{x^{'}})=
-\frac{x^{'}}{x_{0}^{2}}(\frac{x^{'2}}{x_{0}^{2}}-1)<0,
\label{b11}
\end{eqnarray}
since $x_{0}<x^{'} <0$.
The outward normal  vector
is $\overrightarrow{n_{o}}(\overrightarrow{x})=(-1,-\frac{y_{0}}{x_{0}})^{T}$.
The inner product of the the outward normal  vector
 $\overrightarrow{n_{o}}(\overrightarrow{x})$ and the tangent vector $\overrightarrow{f}(\overrightarrow{x^{'}})$
is
\begin{eqnarray}
\overrightarrow{n_{o}}(\overrightarrow{x}).\overrightarrow{f}(\overrightarrow{x^{'}})=
\frac{x^{'}}{x_{0}^{2}}(\frac{x^{'2}}{x_{0}^{2}}-1)<0,
\label{b12}
\end{eqnarray}
since $x^{'} < x_{0} <0$. By Theorem 1, the periodic solution is attractive inward and outward in the third quadrant.
Similarly the attractiveness in the other quadrants can be derived.
In the second quadrant, the inward normal  vector
is $\overrightarrow{n_{i}}(\overrightarrow{x})=(1,\frac{y_{0}}{x_{0}})^{T}$.
The inner product of the the inward normal  vector
 $\overrightarrow{n_{i}}(\overrightarrow{x})$ and the tangent vector $\overrightarrow{f}(\overrightarrow{x^{'}})$
is
\begin{eqnarray}
\overrightarrow{n_{i}}(\overrightarrow{x}).\overrightarrow{f}(\overrightarrow{x^{'}})=
-\frac{x^{'}}{x_{0}^{2}}(\frac{x^{'2}}{x_{0}^{2}}-1)<0,
\label{b13}
\end{eqnarray}
since $x_{0} < x^{'} < 0$.
The inner product of the the outward normal  vector
 $\overrightarrow{n_{o}}(\overrightarrow{x})$ and the tangent vector $\overrightarrow{f}(\overrightarrow{x^{'}})$
is
\begin{eqnarray}
\overrightarrow{n_{o}}(\overrightarrow{x}).\overrightarrow{f}(\overrightarrow{x^{'}})=
\frac{x^{'}}{x_{0}^{2}}(\frac{x^{'2}}{x_{0}^{2}}-1)<0,
\label{b14}
\end{eqnarray}
since $x^{'}< x_{0} <0$.
In the first quadrant, the inward normal  vector
is $\overrightarrow{n_{i}}(\overrightarrow{x})=(-1,-\frac{y_{0}}{x_{0}})^{T}$.
The inner product of the the inward normal  vector
 $\overrightarrow{n_{i}}(\overrightarrow{x})$ and the tangent vector $\overrightarrow{f}(\overrightarrow{x^{'}})$
is
\begin{eqnarray}
\overrightarrow{n_{i}}(\overrightarrow{x}).\overrightarrow{f}(\overrightarrow{x^{'}})=
\frac{x^{'}}{x_{0}^{2}}(\frac{x^{'2}}{x_{0}^{2}}-1)<0,
\label{b15}
\end{eqnarray}
since $x_{0} > x^{'} > 0$.
The inner product of the the outward normal  vector
 $\overrightarrow{n_{o}}(\overrightarrow{x})$ and the tangent vector $\overrightarrow{f}(\overrightarrow{x^{'}})$
is
\begin{eqnarray}
\overrightarrow{n_{o}}(\overrightarrow{x}).\overrightarrow{f}(\overrightarrow{x^{'}})=
-\frac{x^{'}}{x_{0}^{2}}(\frac{x^{'2}}{x_{0}^{2}}-1)<0,
\label{b16}
\end{eqnarray}
since $x^{'}> x_{0} >0$.
In the fourth quadrant, the inward normal  vector
is $\overrightarrow{n_{i}}(\overrightarrow{x})=(-1,-\frac{y_{0}}{x_{0}})^{T}$.
The inner product of the the inward normal  vector
 $\overrightarrow{n_{i}}(\overrightarrow{x})$ and the tangent vector $\overrightarrow{f}(\overrightarrow{x^{'}})$
is
\begin{eqnarray}
\overrightarrow{n_{i}}(\overrightarrow{x}).\overrightarrow{f}(\overrightarrow{x^{'}})=
\frac{x^{'}}{x_{0}^{2}}(\frac{x^{'2}}{x_{0}^{2}}-1)<0,
\label{b17}
\end{eqnarray}
since $x_{0} > x^{'} > 0$.
The inner product of the the outward normal  vector
 $\overrightarrow{n_{o}}(\overrightarrow{x})$ and the tangent vector $\overrightarrow{f}(\overrightarrow{x^{'}})$
is
\begin{eqnarray}
\overrightarrow{n_{o}}(\overrightarrow{x}).\overrightarrow{f}(\overrightarrow{x^{'}})=
-\frac{x^{'}}{x_{0}^{2}}(\frac{x^{'2}}{x_{0}^{2}}-1)<0,
\label{b18}
\end{eqnarray}
since $x^{'}> x_{0} >0$. By Theorem 1, the periodic solution is attractive inward and outward in the
full plane.

\subsection {Attractiveness of stable and unstable manifolds of the equilibria}

For the canonical system
 \begin{eqnarray}
 \left\{\begin{array}{l}
\dot{x}=x,\\
\dot{y}=-y,
\end{array}\right.\label{b19}
\end{eqnarray}
obviously the zero is it's saddle equilibrium,  $x=0$ and $y=0$ are respectively the zero's stable and unstable manifolds.
These manifolds are also the IMs,  $x=0$ is repulsive and $y=0$ is attractive.
Now let's verify their attractiveness and repulsiveness  by Theorem 1.

Firstly let's talk about the attractiveness of the right half of the IM $y=0$, where $x>0$.
Let any point $(x_{0}, 0)$ in the IM $y=0$, where $x_{0}>0$, the tangent passing $(x_{0}, 0)$
is IM $y=0$ itself, thus the normal is $x \equiv x_{0}$. The tangent vector passing any point $(x_{0}, y^{'})$ in the normal,
where $y^{'}>0$,
is $\overrightarrow{f}(\overrightarrow{x^{'}})=(x_{0}, -y^{'})^{T}$.
 The upward normal  vector
is $\overrightarrow{n_{u}}(\overrightarrow{x})=(0, 1)^{T}$.
The inner product of the the upward normal  vector
 $\overrightarrow{n_{u}}(\overrightarrow{x})$ and the tangent vector $\overrightarrow{f}(\overrightarrow{x^{'}})$
is
\begin{eqnarray}
\overrightarrow{n_{u}}(\overrightarrow{x}).\overrightarrow{f}(\overrightarrow{x^{'}})=
-y^{'}<0,
\label{b20}
\end{eqnarray}
since $y^{'} > 0$.
The downward normal  vector
is $\overrightarrow{n_{d}}(\overrightarrow{x})=(0, -1)^{T}$.
The inner product of the the downward normal  vector
 $\overrightarrow{n_{d}}(\overrightarrow{x})$ and the tangent vector $\overrightarrow{f}(\overrightarrow{x^{'}})$
is
\begin{eqnarray}
\overrightarrow{n_{u}}(\overrightarrow{x}).\overrightarrow{f}(\overrightarrow{x^{'}})=
y^{'}<0,
\label{b21}
\end{eqnarray}
since $y^{'} < 0$. Thus the right half of the unstable manifold  $y=0$ is attractive by Theorem 1.
The attractiveness of the left half of the unstable manifold  $y=0$ is can be deduced similarly  by Theorem 1.
So the unstable manifold  $y=0$ is attractive.

Then let's consider the repulsiveness of the upper half of the IM $x=0$, where $y>0$.
Let any point $(0, y_{0})$ in the IM $x=0$, where $y_{0}>0$, the tangent passing $(0, y_{0})$
is IM $x=0$ itself, thus the normal is $y \equiv y_{0}$. The tangent vector passing any point $(x^{'}, y_{0})$ in the normal,
where $x^{'}>0$,
is $\overrightarrow{f}(\overrightarrow{x^{'}})=(x^{'}, -y_{0})^{T}$.
 The rightward normal  vector
is $\overrightarrow{n_{r}}(\overrightarrow{x})=(1, 0)^{T}$.
The inner product of the the rightward normal  vector
 $\overrightarrow{n_{r}}(\overrightarrow{x})$ and the tangent vector $\overrightarrow{f}(\overrightarrow{x^{'}})$
is
\begin{eqnarray}
\overrightarrow{n_{r}}(\overrightarrow{x}).\overrightarrow{f}(\overrightarrow{x^{'}})=
x^{'}>0,
\label{b22}
\end{eqnarray}
since $x^{'}>0$.
The leftward normal  vector
is $\overrightarrow{n_{l}}(\overrightarrow{x})=(-1, 0)^{T}$.
The inner product of the the leftward normal  vector
 $\overrightarrow{n_{l}}(\overrightarrow{x})$ and the tangent vector $\overrightarrow{f}(\overrightarrow{x^{'}})$
is
\begin{eqnarray}
\overrightarrow{n_{l}}(\overrightarrow{x}).\overrightarrow{f}(\overrightarrow{x^{'}})=
-x^{'}>0,
\label{b23}
\end{eqnarray}
since $x^{'} < 0$. Thus the upper half of the stable manifold  $x=0$ is repulsive by Theorem 1.
The repulsiveness of the lower half of the stable manifold  $x=0$ is can be deduced similarly  by Theorem 1.
So the stable manifold  $x=0$ is repulsive.

\subsection {Attractiveness of the "real" IM}

There are the IMs satisfy the definition of IM other than
the equilibria, limit cycle, stable and unstable manifold.
They are called as the real IM.

For the system
 \begin{eqnarray}
 \left\{\begin{array}{l}
\dot{x}=xy^{3},\\
\dot{y}=-y-x-xy^{3},
\end{array}\right.\label{b24}
\end{eqnarray}
obviously $y=-x$ is the IM of system  $(\ref{b24})$.
Let any point $(x_{0},y_{0})$ in the IM, where $y_{0}=-x_{0}$,
the upward normal passing point $(x_{0},y_{0})$ is $y=x-x_{0}+y_{0}=x-2x_{0}$.
Let any point $(x^{'},y^{'})$ in the upper normal, where $y^{'}=x^{'}-2x_{0}$,
the tangent vector passing
point $(x^{'},y^{'})$ is $\overrightarrow{f}(\overrightarrow{x^{'}})=(x^{'}y^{'3},-y^{'}-x^{'}-x^{'}y^{'3})^{T}$,
where $x^{'}>x_{0}$. The upward normal vector is $\overrightarrow{n_{u}}(\overrightarrow{x})=(1, 1)^{T}$.
The inner product of the the upward normal  vector
 $\overrightarrow{n_{u}}(\overrightarrow{x})$ and the tangent vector $\overrightarrow{f}(\overrightarrow{x^{'}})$
is
\begin{eqnarray}
\overrightarrow{n_{u}}(\overrightarrow{x}).\overrightarrow{f}(\overrightarrow{x^{'}})=
-x^{'}-y^{'}=2 (x_{0}-x^{'})<0,
\label{b25}
\end{eqnarray}
since $x^{'}>x_{0}$.
The downward normal vector is $\overrightarrow{n_{d}}(\overrightarrow{x})=(-1, -1)^{T}$.
The inner product of the the downward normal  vector
 $\overrightarrow{n_{d}}(\overrightarrow{x})$ and the tangent vector $\overrightarrow{f}(\overrightarrow{x^{'}})$
is
\begin{eqnarray}
\overrightarrow{n_{d}}(\overrightarrow{x}).\overrightarrow{f}(\overrightarrow{x^{'}})=
x^{'}+y^{'}=-2(x_{0}-x^{'})<0,
\label{b25}
\end{eqnarray}
since $x_{0}>x^{'}$.
Thus the IM  $y=-x$ is attractive upward and downward by Theorem 1.

\section {Conclusion and Discussion}

In this paper, we first present an operable but not abstract, universal  and simple theory,
Theorem 1, on the attractiveness of the
IMs, i.e., for any point $\overrightarrow{x}$ in the IM
$M$, $\overrightarrow{n}(\overrightarrow{x})$ is the normal passing by $\overrightarrow{x}$, and $\forall \overrightarrow{x^{'}} \in  \overrightarrow{n}(\overrightarrow{x})$, if the tangent $\overrightarrow{f}(\overrightarrow{x^{'}})$ of the orbits of the
dynamical system intersects at  obtuse angle with the normal  $\overrightarrow{x}$, i.e., $\overrightarrow{f}(\overrightarrow{x^{'}}).
\overrightarrow{n}(\overrightarrow{x}) < 0$, then the IM $M$ is attractive. The
conclusion of the repulsiveness  is also obtained similarly. Some illustrative
examples of the IMs, such as equilibria, periodic solution, stable and unstable manifolds,
other IM are presented to support our result. This method
is simple and universal for different kinds of IMs.

This method provides an operable theory to
consider the attractiveness of the
IMs.  But there are two problems to be solved.
If the considered systems or the IMs are higher-dimensional,
including the delayed differential dynamical systems,
the    normal  vector
 $\overrightarrow{n_{d}}(\overrightarrow{x})$ and the tangent vector $\overrightarrow{f}(\overrightarrow{x^{'}})$
 will be too difficult to be obtained. Their inner product is
 hard to be derived. On the other hand, this method can be applied only in the case that
 the expression of IM is analytically known.  But the closed form of IM, such as the synchronization
 manifolds in the generalized synchronization,  is usually
 hard to be solved since it is in fact the solution of a first-ordered quasi-linear PDEs.
 But we can obtain their approximated expression and consider
 their attractiveness or repulsiveness by Theorem 1.
These  will my next work.

\section*{Acknowledgment}

The author would like to acknowledge the financial support for this research
 via the Natural National Science   Foundation of China (No. 10702065). He
 also thanks the reviewers for their valuable reviews and
 suggestions.






\begin{thebibliography}{1999}



\bibitem{ockendon}  J. R. Ockendon, S. D. Howison, A. A. Lacery, et al, Applied Partial Equations, Oxford: Oxford University Press  (2003).

\bibitem{fenichel} N. Fenichel,   Persistence and smoothness of invariant manifolds for flows, Indiana University Mathematics Journal, 21(3) (1971) 193¨C225.

 \bibitem{bronstein}
I. U. Bronstein, A. Y. Kopanskii, Smooth invariant manifolds and normal forms, Singapore: World Scientific (1994).

\bibitem{wiggins}
 S. Wiggins, Normally Hyperbolic Invariant Manifolds in Dynamical Systems, Applied Mathematical Sciences (1994).

\bibitem{chow} S. N. Chow, W. Liu,  Synchronization, stability and normal hyper-bolicity, Resenhas IME-USP 3(1997) 139¨C158.

\bibitem{josic} K. Josi$\acute{c}$,  Synchronization of chaotic systems and invariant manifolds,  Nonlinearity,  13(4)(2000)  1321-1336.

 \bibitem{gorban} A. N. Gorban, I. V. Karlin, A. Y. Zinovyev, Constructive methods of invariant manifolds for kinetic problems,  Physics Reports   396(4-6)(2004)  197-403.

\bibitem{lyapunov} A. M. Lyapunov, General Problem of the Stability of Motion, Taylor and  Francis Books Ltd  (1992).
\end{thebibliography}




\end{document}